\documentclass{amsart}
\usepackage[colorlinks,pdfstartview={XYZ null null null}]{hyperref}
\usepackage{url}

\newcommand{\Z}{\mathcal{Z}}            
\newcommand{\W}{\mathcal{W}}            
\newcommand{\A}{\mathcal{A}}            
\newcommand{\N}{\mathcal{N}}            
\newcommand{\E}{\mathcal{E}}            
\newcommand{\ID}{\mathbb{D}}            
\newcommand{\IT}{\mathbb{T}}            
\newcommand{\IC}{\mathbb{C}}            
\DeclareMathOperator{\diam}{diam}
\DeclareMathOperator{\re}{Re}
\let\union\cup              

\let\Union\bigcup

\newcommand*{\parenth}[1]{\textup{(}#1\textup{)}}

\newcommand{\defeq}{\stackrel{\mathrm{def}}{=}}
\newcommand{\st}{:}
\let\eps\epsilon
\let\phi\varphi
\let\z\zeta

\let\term\textit

\newtheorem{theorem}{Theorem}[section]

\newtheorem{Athm}{Theorem}

\newtheorem{fact}{Fact}

\theoremstyle{remark}

\numberwithin{equation}{section}

\renewcommand{\labelenumi}{{\upshape(\theenumi)}}

\begin{document}

\title[Interpolating Sequences]%
{Finite Unions of Interpolating Sequences for Hardy Spaces}

\author{Daniel H. Luecking}
\address{Department of Mathematical Sciences\\
         University of Arkansas\\
         Fayetteville, Arkansas 72701}
\email{luecking@uark.edu}
\date{June 16, 2005}

\begin{abstract}
  A sequence which is a finite union of interpolating sequences for
  $H^\infty$ have turned out to be especially important in the study of
  Bergman spaces. The Blaschke products $B(z)$ with such zero sequences
  have been shown to be exactly those such that the multiplication $f
  \mapsto fB$ defines an operator with closed range on the Bergman
  space. Similarly, they are exactly those Blaschke products that
  boundedly divide functions in the Bergman space which vanish on their
  zero sequence. There are several characterizations of these sequences,
  and here we add two more to those already known. We also provide a
  particularly simple new proof of one of the known characterizations.
  One of the new characterizations is that they are interpolating
  sequences for a more general interpolation problem.
\end{abstract}

\maketitle

\section{Introduction}\label{introduction}

Let $\ID$ denote the open unit disk and $\IT$ its boundary. Let $m$
denote the normalized arc length measure on $\IT$. For $0< p \le \infty$
let $L^p$ denote the usual Lebesgue space of all measurable functions
$f$ on $\IT$ for which $|f|^p$ is integrable with respect to $dm$, and
let $H^p$ denote the usual Hardy space consisting of analytic functions
$f$ on $\ID$ such that $\sup_{0<r<1} M_p(f,r) < \infty$, where
\begin{equation*}
  M_p(f,r) =
  \begin{cases}
  \left( \frac{1} {2\pi} \int_0^{2\pi} |f(re^{i\theta})|^p
    \,d\theta \right)^{1/p}& 0 < p < \infty\\
    \sup_{0 \le \theta < 2\pi} |f(re^{i\theta})| & p = \infty
  \end{cases}
\end{equation*}
We denote both the norm on $L^p$ and the norm on $H^p$ (defined by the
supremum above) by $\| f \|_{H^p}$. It is well known that functions in $H^p$
have radial limits $\lim_{r\to 1} f(re^{i\theta})$ for $m$-almost all
$\theta$ and that this defines an isometry between $H^p$ and a closed
subspace of $L^p$. We abuse notation by writing $f$ for both the
function $f(z)$ on $\ID$ and its limit $f(e^{i\theta})$ on $\IT$. We
also abuse the terminology by calliing $\| \cdot \|_{H^p}$ a norm even in
the cases $0 < p < 1$.

Let $\psi(z,\zeta)$ denote the \term{pseudohyperbolic metric}:
\begin{equation*}
  \psi(z,\zeta) = \left|\frac{z-\zeta}{1-\bar\zeta z}\right|.
\end{equation*}
We will use $D(z,r)$ for the \term{pseudohyperbolic disk} of radius $r$
centered at $z$, that is, the ball of radius $r<1$ in the
pseudohyperbolic metric.

Let $\Z= \{ z_k : k=1,2,3,\dots \}$ be a sequence in $\ID$ without limit
points in $\ID$. Define the space of sequences $l^p_\Z$, $0 < p \le
\infty$, to be all those $w = (w_k)$ such that $\| w \|_{l^p_\Z} <
\infty$, where
\begin{equation*}
  \| w \|_{l^p_\Z}^p =
    \begin{cases}
      \sum_j |w_k|^p(1 - |z_k|^2) & \text{if $p < \infty$}\\
      \sup_j |w_k|                & \text{if $p = \infty$}
    \end{cases}
\end{equation*}

The usual $H^p$ interpolation problem for $\Z$ is the following: given a
sequence $w = (w_k) \in l^p_\Z$ find a function $f\in H^p$ such that
$f(z_k) = w_k$ for all $k$. The sequence $\Z$ is called an
\term{interpolating sequence for $H^p$} if every such interpolation
problem has a solution. Interpolating sequences for $H^\infty$ were
characterized by L.~Carleson in \cite{Car58} and this was extended to
all other $H^p$ by H.~S.~Shapiro and A.~Shields in \cite{SS61}. Let
$T_p$ denote the interpolation operator: $T_p(f)$ is the sequence
$(f(z_j) \st j\ge 1)$. Note that the definition of interpolating
sequence requires $l^p_\Z \subset T_p(H^p)$ but does not required
$T_p(H^p) \subset l^p_\Z$, however in \cite{SS61} it was shown that the
former implies the latter.

The characterization of interpolating sequences is the same for all
$H^p$: a sequence $\Z$ is interpolating for $H^p$ ($0 < p \le \infty$)
if and only if $\Z$ is \term{uniformly separated}, which means there
exists $\delta > 0$ such that
\begin{equation*}
  \inf_k \prod_{j \ne k} \left| \frac{z_k - z_j}{1 - \bar z_j
    z_k} \right| > \delta.
\end{equation*}
This clearly implies that $\Z$ is \term{uniformly discrete}, which is to
say $\psi(z_j,z_k) \ge \delta$ for all $j \ne k$. Closely associated
with $\psi$ are the \term{Moebius transformations}
\begin{equation*}
    \phi_\zeta(z) = \frac{\zeta-z}{1-\bar\zeta z},\quad \zeta \ne 0.
\end{equation*}
Clearly $\psi(z,\zeta) = |\phi_\zeta(z)|$.

Since an interpolating sequence for $H^p$ is a zero sequence for some $f
\in H^p$, it follows that an interpolating sequence $\Z = \{ z_1, z_2,
\dots \}$ is a \term{Blaschke sequence}, that is $\sum 1 -|z_k|^2 <
\infty$. The \term{Blaschke product} $B_\Z$ determined by $\Z$ is
\begin{equation*}
    B_\Z(z) = z^m\prod_k \frac{\bar z_k}{|z_k|}\phi_{z_k}(z)
\end{equation*}
where $m$ is the number of repetitions of $0$ in the sequence $\Z$. This
product converges uniformly on compact sets in $\ID$ if and only if $\Z$
is a Blaschke sequence. It provides a bounded analutic function that
vanishes precisely on $\Z$, with the order of each zero equal to the
number of repetitions of that point in the sequence $\Z$. When we say
``$\Z$ is the zero sequence of $f$'' we will always mean that the order
of the zero of $f$ at $\zeta$ equals the repetition of $\zeta$ in the
sequence $\Z$. Also, ``$f$ vanishes on $\Z$'' will mean $f$
vanishes at $\zeta$ to order at least the repetition of $\zeta$ in $\Z$.
Similar interpretations are to be given to statements such as ``$\Z
\subset \W$'' when $\Z$ and $\W$ are sequences in $\ID$.

The \term{Bergman space} $A^p$ is defined to be the the collection of
analytic functions $f$ on $\ID$ that belong to $L^p(dA)$ where $dA$
denotes area measure on $\ID$. The norm on $A^p$ will be denoted $\|
\cdot \|_{A^p}$. For $\alpha > -1$ let $dA_\alpha$ denote the measure
defined by $dA_\alpha(z) = (1 + |z|^2)^\alpha\,dA(z)$. Then the
\term{weighted Bergman space $A^{p,\alpha}$} is the set of all analytic
functions in $\ID$ such that
\begin{equation*}
    \| f \|_{A^{p,\alpha}} \defeq \left( \int_{\ID} |f|^p \,dA_\alpha
        \right)^{1/p} < \infty.
\end{equation*}

In a recent paper, P.~Duren and A~Schuster \cite{DS02} proved the
following equivalence of five conditions. Some terms used here will be
defined only in the next section. The unit point mass at a point a is
denoted by $\delta_a$.

\begin{Athm}\label{D-S}
\renewcommand\theenumi{\roman{enumi}}
\renewcommand\labelenumi{\parenth\theenumi}
  For a Blaschke sequence $\Z = \{ z_k \}$ of points in $\ID$, the
  following five statements are equivalent.
  \begin{enumerate}
    \item $\Z$ is a finite union of uniformly separated sequences.
        \label{FUS}
    \item $\sum_{k=1}^{\infty}(1 - |z_k|^2)\delta_{z_k}$ is a Carleson
        measure.
    \item $\sup_{z\in\ID} \sum_{k=1}^{\infty}(1 - |\phi_{z_k}(z)|^2) <
        \infty$. \label{UBC}
    \item The associated Blaschke product $B$ is a universal divisor of
        $A^p$; that is, $f/B \in A^p$ for every function $f \in A^p$
        that vanishes on $\Z$. \label{UD}
    \item The operator $M_B$ of multiplication by the associated Blaschke
        product $B$ is bounded below on $A^p$; that is, there is a
        constant $c > 0$ such that $\| Bf \|_{A^p} \ge c\| f \|_{A^p}$
        for every function $f \in A^p$. \label{BBM}
  \end{enumerate}
\end{Athm}

In this paper we will provide a simpler proof of the fact that
(\ref{UBC}) implies (\ref{UD}). In addition, we prove the following
sixth and seventh equivalent properties. The term \term{general
interpolating sequence} will be defined later. Suffice it to say for now
that it is a sequence for which a particular interpolation problem
always has a solution in the indicated space of functions.

\begin{theorem}\label{main}
\renewcommand\theenumi{\roman{enumi}}
\renewcommand\labelenumi{\parenth\theenumi}
  For a Blaschke sequence $\Z = \{ z_k \}$ the following are equivalent.
  \begin{enumerate}
    \item $\Z$ is a finite union of uniformly separated sequences.
          \addtocounter{enumi}{4}
    \item For the associated Blaschke product $B$, the zero function
        does not belong to the closure of $\{ B\circ \phi_z \mid z\in
        \ID \}$ in the topology of uniform convergence on compacta.
        \label{UNZ}
    \item $\Z$ is a general interpolating sequence for $H^p$. \label{GIS}
  \end{enumerate}
\end{theorem}

In statements (\ref{UD}), (\ref{BBM}) and (\ref{GIS}) the value of $p$
is ambiguous. As in Duren and Schuster's paper \cite{DS02}, it will
follow from the proof that if any of the statements hold for one value of
$p$ then they all hold for all values of $p$.

We will prove that statement (\ref{BBM}) of the Theorem~\ref{D-S}
implies (\ref{UNZ}) of Theorem~\ref{main} and that this implies
(\ref{FUS}). And then we will prove that (\ref{GIS}) of
Theorem~\ref{main} is equivalent to the rest of the conditions.

\section{Background}

A measure $\mu\ge 0$ on $\ID$ is called a \term{Carleson measure} if
there is a finite constant $C \ge 0$ such that for every arc $I \subset
\IT$
\begin{equation*}
  \mu(S_I) \le Cm(I)
\end{equation*}
where $S_I = \{ z\in\ID \st z/|z| \in I \text{ and } 1 - |z| < m(I) \}$
is the \term{Carleson square based on $I$}. The infimum of the set of
$C$ for which this inequality holds is denoted $\| \mu \|_*$. If any $0
< p < \infty$ is given, Carleson measures are precisely those for which
there exists a constant $C$ with
\begin{equation*}
  \int |f|^p \,d\mu \le C\| f \|_{H^p}^p \quad \text{for all $f\in H^p$}.
\end{equation*}
In fact, there exist absolute constants $a$ and $B$ such that $a\| \mu
\|_* \le C \le  B\| \mu \|_*$.

Given the sequence $\Z = \{ z_j \st j \ge 1 \}$ let $\mu_\Z$ denote the
measure $\sum (1 - |z_j|^2)\delta_{z_j}$. If $\Z$ is an interpolating
sequence for $H^p$ then, because it is a Blaschke sequence, $\mu_\Z$ is
a finite measure. Another characterization of interpolating sequences
(essentially also obtained in \cite{Car58} and \cite{SS61}) is that they
must be uniformly discrete and $\mu_\Z$ must be a Carleson measure. A
sequence where $\mu_\Z$ is a Carleson measure will be called a
\term{Carleson sequence}. It is well know (see \cite{McDS79}) that they
can be written as a finite union of uniformly discrete sequences, and so
they are just the finite unions of interpolating sequences.

We can restate the definition of interpolating sequence in terms of the
Banach spaces $H^p$ and $l^p_\Z$. An interpolating sequences is one for
which $l^p_\Z$ is contained in the range of the interpolation map $T_p:
H^p \to l^p_\Z$. As already remarked, although it is not immediately
obvious, once this range contains $l^p_\Z$, it is necessarily contained
in it and the closed graph theorem implies $T_p$ is necessarily
continuous. Then, by the open mapping theorem, if $\Z$ is an
interpolating sequence for $H^p$, there is a finite constant $M > 0$
such that every $w\in l^p_\Z$ is interpolated by some $f_w \in H^p$ with
$\| f_w \|_\infty \le M\| w \|_\infty$. This fact can actually be
obtained without the open mapping theorem, using a normal families
argument. The smallest such constant $M$ is called the
\term{interpolation constant} of the sequence $\Z$.

An interpolating sequence necessarily has no repeated points and must be
a Blaschke sequence. A Blaschke product whose zero sequence is an
interpolating sequence is called an \term{interpolating Blaschke
product}.

\section{Finite Products of interpolating Blaschke}\label{fpibp}

The earliest mention of finite products of interpolating Blaschke
products that I am aware of is in the paper of McDonald and Sundberg
\cite{McDS79}, where they were part of the invertibility criterion for
Toeplitz operators with analytic symbol on the Bergman space of $\ID$.

In \cite{DS02}, P.~Duren and A.~Schuster compiled most of the known
characterizations of these Blaschke products (equivalently, of Carleson
sequences). They also provided a unified proof of these characterization
in the sense that each condition was shown to imply the next in (almost)
a cycle, and the proof of each implication was relatively painless.

In this section we discuss those conditions, plus one of the conditions
in Theorem~\ref{main}. Let $\Z$ be a Blaschke sequence (of distinct
points) and let $B=B_\Z$ be the associated Blaschke product. For each $z_j
\in \Z$, let $B_j(z)$ be the Blaschke product associated with
$\Z\setminus\{z_j\}$. That is
\begin{equation*}
    B_j(z) = \prod_{k \ne j} \frac{|z_k|}{z_k}\frac{z_k-z}{1-\bar z_k z}
\end{equation*}
Then a Blaschke sequence $\Z$ is uniformly separated if and only if
\begin{equation}\label{US}
  \eps = \inf_{j} |B_j(z_j)| > 0.
\end{equation}
The following is easily seen to be equivalent to \eqref{US}:
\begin{equation*}
    \eps = \inf_{j} (1 - |z_j|^2)|B'(z_j)| > 0.
\end{equation*}
This follows from the product rule for derivatives and the formula for the
derivative of $\phi_{z_j}$.

While the composition of a Blaschke product with a conformal selfmap of
$\ID$ is again a Blaschke product or a constant multiple of one, the
Blaschke condition on the zero sequence is not uniform.
If $B$ is the Blaschke product with zero sequence $\Z$, and $\phi$ is a
conformal selfmap, then $\phi(\Z)$ is the zero sequence of $B\circ
\phi^{-1}$. The Blaschke condition of $\phi(\Z)$ is $\sum_{j} (1 -
|\phi(z_j)|^2) < \infty$. We will say that $\Z$ satisfies a \term{uniform
Blaschke condition} if
\begin{equation*}
  \sup_\phi \sum_{j} (1 - |\phi(z_j)|^2) < \infty
\end{equation*}
where the supremum is taken over all conformal maps of $\ID$ onto
itself.

A Blaschke product $B$ with zero sequence $\Z$ is called a
\term{universal divisor} for $A^p$ if there is a constant $C$ such that
$\| f/B \|_{A^p} \le C\| f \|_{A^p}$ whenever $f\in A^p$ vanishes on $\Z$
(counting multiplicity). By the closed graph theorem, the existence of
the constant is implied by the seemingly weaker statement that $f/B \in
A^p$ whenever $f\in A^p$ vanishes on $\Z$.

We define the operator $M_B$ on $A^p$ by $M_B f = Bf$. It clearly
satisfies $\|M_B f\|_{A^p} \le \| f \|_{A^p} $ for any
$0<p<\infty$. We will say that $M_B$ is \term{bounded
below on $A^{p}$} if there exists a constant $c > 0$ such that
$\| M_B f \|_{A^p} \ge c\| f \|_{A^p}$.

It is quite easy to create Blaschke products $B$ such that for some
sequence $\phi_{\zeta_n}$ of conformal selfmaps of $\ID$ one has $B\circ
\phi_{\zeta_n}$ tending to $0$ uniformly on compact subsets of $\ID$.
For example, let $\{ \zeta_n\}$ be a sequence of distinct points with
$|\zeta_n| \to 1$ and satisfying $\sum n(1-|\zeta_n|) < \infty$. Let
$\Z$ be the sequence containing each $\zeta_n$ repeated $n$ times. Then
if $B$ is the associated Blaschke product, we see that
$|B(\phi_{\zeta_n}(z))|$ has a zero at $0$ of order $n$, and so it has
the form $z^n g(z)$ with $\| g \|_\infty =1$. Therefore, it tends to $0$
uniformly on compact. (This example also has the property the
$\sum_{a\in \Z} (1 - |\phi_{\zeta_n}(a)|^2) \ge n$, so it idoes not
satisfy the uniform Blaschke condition.)

We will say that $B$ is \term{uniformly bounded away from zero} (or
\term{uniformly nonzero} for short) if the above described phenomenon
does not occur. That is, for any sequence $\phi_n(z)$ of conformal
selfmaps of $\ID$ the compositions $B\circ\phi_n$ do not tend to 0
uniformly on compacta.

Now we can state part of the content of Theorems~\ref{D-S} and
\ref{main} as the equivalence of the following conditions on $\Z$ and
$B=B_\Z$:
\begin{enumerate}
  \item $\Z$ is a finite union of uniformly separated sequences.
  \item $\Z$ is a Carleson sequence.
  \item $\Z$ satisfies the uniform Blaschke condition.
  \item $B$ is a universal divisor for $A^p$.
  \item $M_B$ is bounded below on $A^p$
  \item $B$ is uniformly nonzero.
\end{enumerate}

In the next section we will provide a different proof that the uniform
Blaschke condition on $\Z$ implies that $B_\Z$ is a universal divisor
for $A^p$. In fact the proof works (and is no less simple) for weighted
Bergman spaces $A^{p,\alpha}$. The proof that $M_B$ is bounded below for
any universal divisor is trivial in any space. The fact that $B$ is
uniformly nonzero when $M_B$ is bounded below also has the same proof
for any weighted Bergman space, and so the equivalence of all the
conditions is not only independent of $p$ but it is also true for any
weighted $A^{p,\alpha}$ and is independent of $\alpha$.

The last condition is probably the minimum conceivable necessary
condition.  All the other conditions (and the constants implicitly
involved in them) are easily shown to be conformally invariant.
Moreover, if we take each condition and fix the constant, it clearly
implies that no limit of such Blaschke products can be zero. The proof
will show that we get a slightly stronger result in one direction.
The last condition can be replaced by the following: {\it For any
subsequence $\zeta_n$ of $\Z$ the compositions $B\circ\phi_{\zeta_n}$ do
not tend to zero uniformly on compacta}. In other words, for
sufficiency, we do not need to test all sequences of Moebius
transformations, but only those of the form $\phi_\zeta$ for $\zeta \in
\Z$.

\section{Proofs}\label{proof}

We start with the simplified proof that (\ref{UBC}) implies (\ref{UD})
in Theorem~\ref{D-S}. It is simpler in that it doesn't need any of the
results of C.~Horowitz' paper \cite{Hor74}, using only basic
inequalities.

Let $f$ be a function in $A^p$ vanishing on $\Z$ and suppose $\Z$
satisfies (\ref{UBC}). Let $\W=\{ \zeta_1, \zeta_2,\dots \}$ be the
zero sequence of $f$. Assume temporarily that $f(0) \ne 0$. We apply
Jensen's formula writing it in the form:
\begin{equation*}
    \log |f(0)| + \sum_{j} \log \frac{r}{|\zeta_j|} \chi_{[|\zeta_j|,1)}(r) =
        \frac{1}{2\pi}\int_0^{2\pi} \log |f(re^{i\theta})| \,d\theta.
\end{equation*}
If we sum only over the points in $\Z$ this equality becomes an
inequality:
\begin{equation*}
    \log |f(0)| + \sum_{j} \log \frac{r}{|z_j|} \chi_{[|z_j|,1)}(r) \le
        \frac{1}{2\pi}\int_0^{2\pi} \log |f(re^{i\theta})| \,d\theta.
\end{equation*}
Now multiply this inequality by $2r\,dr$ and integrate with respect to
$r$ from $0$ to $1$:
\begin{equation*}
    \log |f(0)| + \sum_{j} \left[ \log \frac{1}{|z_j|} - \frac{1 -
    |z_j|^2}{2} \right] \le
        \frac{1}{\pi}\int_{\ID} \log |f(w)| \,dA(w).
\end{equation*}
Now fix $\zeta\in\ID$ with $f(\zeta) \ne 0$ and apply the above to
$f\circ\phi_\z$
where, as before, $\phi_\z$ is the conformal map that interchanges $0$ and
$z$. The mapping $\phi_\z$ is it's own inverse so the zeros of
$f\circ\phi_\z$ are $\phi_\z(z_j)$:
\begin{equation*}
    \log |f(\z)| + \sum_{j} \left[ \log \frac{1}{|\phi_\z(z_j)|} -
    \frac{1 - |\phi_\z(z_j)|^2}
        {2} \right] \le \frac{1}{\pi}\int_{\ID} \log |f(\phi_\z(w))|
        \,dA(w).
\end{equation*}
Now we  make a change of variable $w \mapsto \phi_\z(w)$ in the
integral on the right, and rearrange the parts of the inequality.
\begin{equation*}
    \log |f(\z)| + \log \frac{1}{|B(\z)|}
        \le \sum_j\frac{1 - |\phi_\z(z_j)|^2}{2} + \frac{1}{\pi}
            \int_{\ID} \log |f(w)| |\phi_\z'(w)|^2 \,dA(w).
\end{equation*}
Applying condition~(\ref{UBC}), the sum on the right hand side is
bounded by some constant $C$.
Multiply this inequality by $p/2$, apply $\exp$ to both
sides, and use the inequality of the geometric and arithmetic mean on
the right side (the measure $|\phi_\z'(w)|^2\,dA(w)/\pi$ has total
measure $1$):
\begin{equation*}
    \left| \frac {f(\z)}{B(\z)} \right|^{p/2} \le \frac{e^{Cp/2}}{\pi}
        \int_{\ID} |f(w)|^{p/2} |\phi_\z'(w)|^2 \,dA(w).
\end{equation*}
Note that the right hand side can be seen as the application to an $L^2$
function (namely $|f|^{p/2})$ of an integral operator with a familiar
kernel:
\begin{equation*}
    K(\z,w) = |\phi_\z'(w)|^2 = \frac{(1 - |\z|^2)^2}{|1 - \bar\z w|^4}.
\end{equation*}
Condition (\ref{UD}) will follow if we can show that this kernel
$K(\z,w)$ defines a bounded integral operator on $L^2(dA)$. This was done
in \cite{Lue96} using what is sometimes called the \term{Schur method},
sometimes the \term{Forelli-Rudin method}. It involves only H\"older's
inequality and basic estimates of the integrals involved (see
\cite{FR74}). In \cite{Lue96} it was used to prove the boundedness of
this integral operator in a wide variety of spaces. In particular, it is
bounded in the weighted space $L^2(dA_\alpha)$. Thus
condition~(\ref{UD}) of Theorem~\ref{D-S} (for any weighted Bergman
space) follows from condition~(\ref{UBC}).

Now let us prove part of Theorem~\ref{main} by showing that $B$ is
uniformly nonzero (condition~(\ref{UNZ})) whenever $M_B$ is bounded
below on any $A^p$ (condition~\ref{BBM}).

Let $\phi$ be a conformal self-map of $\ID$ and consider the function
$f(z) = \phi'(z)^{2/p}$, which has norm $\| f \|_{A^p}^p = \int |\phi'|^2
\,dA = \int 1\,dA $ by a change of variables. Applying
condition~(\ref{BBM}) we get
\begin{equation*}
    \int |B\circ\phi|^p \,dA = \int |B|^p|\phi'|^2 \,dA =
        \int |Bf|^p \,dA = \| Bf \|_{A^p}^p \ge c^p \| f \|_{A^p} = c^p\pi.
\end{equation*}
If $\phi_n$ is a sequence of conformal maps, the above holds with $\phi$
equal to any one of them. Since $B\circ\phi_n$ are uniformly bounded
by $1$, if it were to converge to $0$ pointwise, the Lebesgue theorem would
imply that their $L^p(dA)$ norms tend to 0, contradicting the above
inequality. Thus $B$ must be uniformly nonzero.

Finally, assuming $B_\Z$ is uniformly nonzero (condition~(\ref{UNZ})) we
wish to prove $\Z$ is a finite union of interpolating sequences
(condition~(\ref{FUS})).

First we show that condition~(\ref{UNZ}) implies we can factor $B$ into
a finite number of Blaschke products $F_j$ such that the zero sequence
$\Z_j$ of $F_j$ satisfies $\psi(a_k,a_n) > 1/2$ for all $a_k, a_n \in
\Z_j$, $k \ne n$.

Toward this end we show that there is an upper bound on the number of
zeros of $B$ contained in the pseudohyperbolic disks $D(\zeta, 1/2)$ as
$\zeta$ varies in $\ID$. Let $N(\zeta)$ denote the number of zeros $a_k$
of $B$ in $D(\zeta,1/2)$. For any sequence $\zeta_n$ in $\ID$, let
$\phi_{\zeta_n}$ be the conformal self-map that interchanges 0 and
$\zeta_n$ and observe that $B\circ\phi_{\zeta_n}$ has $N(\zeta_n)$ zeros
in $|z| < 1/2$. If there is no upper bound on $N(\zeta)$, choose a
sequence $\zeta_n$ with $N(\zeta_n) \to \infty$. Then we have
$B\circ\phi_{\zeta_n} \to 0$ uniformly on compacta. This contradicts
(\ref{UNZ}), so there must exist an $N$ such that $N(\z) \le N$ for all
$\z \in \ID$.

It is well-known that $\Z$ can now be partitioned into finitely many (in
fact $N$) subsequences such that each subsequence has a separation of at
least $1/2$ (in the pseudohyperbolic metric) between points.
Thus
$\Z = \union_1^N \Z_j$ and each $\Z_j$ satisfies $\psi(a,b) > 1/2$ for
all pairs $a\ne b$ in $\Z_j$. This induces a factorization of $B =
\prod_{j=1}^{M} F_j$. Let $F$ be one of the factors and note that $F$
must also be uniformly bounded away from zero. We want to show that
the zero set of $F$ is uniformly separated, so assume it is not. Then
there is a sequence $a_k$ of zeros of $F$ such that $(1 -
|a_k|^2)|F'(a_k)| \to 0$. Without loss of generality we may assume the
sequence $G_k = F\circ\phi_{a_k}$ converges uniformly on compact sets to some
function $G$. In the disk $|z| < 1/2$ each $G_k$ has exactly one zero of
order 1: at the origin. By Hurwitz' theorem, $G$ must either have exactly
one zero of order 1 (at the origin) or be identically zero. The latter
is ruled out because $F$ is uniformly nonzero. Therefore,
$$
    (1 - |a_k|^2)|F'(a_k)| = G_k'(0) \to G'(0) \ne 0,
$$
and we have a contradiction.

\section{General interpolation}

In this section we set up the notions of interpolation we call
\term{general} interpolation. This is nearly identical to the
corresponding notion for Bergman spaces introduced in \cite{Lue04b},
except for the norm imposed on the sequence space. Therefore we will
omit many of the (identical) proofs and most of the justification for
certain choices.

For simplicity of notation, let $(S)_\eps$ denote the
$\eps$-neighborhood of a set $S$, in the pseudohyperbolic metric.
That is $(S)_\eps = \Union_{s\in S} D(s,\eps)$.

Let $\Z$ be any sequence in $\ID$ without limit points and let us
suppose it has been partitioned into finite sets $\Z_k$ with the
property that the pseudohyperbolic diameters of $\Z_k$ are bounded away
from $1$. Let there be given open sets $G_k$ with the following
properties: There is an upper bound $R < 1$ on the diameters of $G_k$
and there is an $\eps > 0$ such that $(\Z_k)_\eps \subset G_k$.
One could in fact take $G_k = (\Z_k)_\eps$, but we wish to start as
general as possible.

Now we define a norm $\| \cdot \|_{G_k}$ for analytic functions on $G_k$
to be the the supremum norm. That is, $\| f \|_{G_k} = \sup_{z\in G_k}
|f(z)|$. This is probably the simplest choice of norms, but quite a lot
of other choices produce the same meaning of general interpolation.
Other norms that work (that is, ensure the validity of succeeding
arguments) are $L^p$ averages on $G_k$ with respect to area or with
respect to arc length on the boundary (assuming the pseudohyperbolic
perimeter of the $G_k$ are bounded). We stick with the notation $\|
\cdot \|_{G_k}$ to emphasize that it need not always be the sup norm.

Let $\A_k$ be the space of all bounded analytic functions $f$ on $G_k$
with the sup norm and let $\N_k$ be the closed subspace of $\A_k$ of
all functions that vanish on $\Z_k$ to order given by the repetition
within $\Z_k$. Let $\E_k$ denote the quotient space $\A_k/\N_k$ with the
quotient norm, denoted by $\| \cdot \|_k$.

Let $d_k$ denote the Euclidean distance from $G_k$ to the boundary of
$\ID$ and for $0<p<\infty$ define $X^p$ to be the space of all sequences
$(w_k)$ with $w_k \in \E_k$ that satisfy $\sum \| w_k \|_k^p d_k <
\infty$. Let $X^\infty$ denote the sequences $(w_k)$ with $\sup \| w_k
\|_k < \infty$. The notation $X^p$, for the sake of simplicity,
obviously suppresses a great many details involved in the above
description. The obvious norm in this space is denoted $\| \cdot \|_{X^p}$

We say a function $f \in H^p$ \term{interpolates} a sequence $(w_k)$ if
for each $k$, $f|_{G_k}$ lies in the equivalence class $w_k$. This is
equivalent to saying that at points of $\Z$, $f$ has the values (and
values of derivatives, etc.) that are defined by $w_k$. Thus, this is
truly interpolation. What makes it \term{general} is the more
general sequence space $X^p$.

For example, if each $\Z_k = \{ z_k \}$ is a singleton, the spaces
$\E_k$ are one dimensional and therefore essentially $\IC$. The
equivalence class of $f \in \A_k$ is uniquely determined by the value of
$f$ at the point $z_k$. It is not hard to show that $\| f + \N_k \|_k$
is equivalent to $|f(z_k)|$, and so the above scheme is exactly the
usual notion of interpolation (evaluating at points). The requirement
that $G_k$ have pseudohyperbolic diameters bounded away from $1$ implies
that $d_k$ is equivalent to $1 - |a_k|$ for any choice of $a_k \in G_k$
so the space $X^p$ is essentially $l^p_\Z$.

For another example, let $n_k$ be a bounded sequence of nonnegative
integers and let each $\Z_k$ be a single point $z_k$ repeated $n_k$
times. Then the equivalence class of $f\in \A_k$ is determined by the
values of $f$, $f'$,\dots,$f^{(n_k-1)}$ at $z_k$ and we get the usual
notion of \term{multiple interpolation} where we obtain our sequence by
evaluating a function and its derivatives (up to some order) at $z_k$.
It is not hard to show that the norm of the equivalence class $f + \N_k$
is essentially $|f(z_k)| + |f'(z_k)|d_k + \cdots + |f^{(n_k - 1)}|
d_k^{n_k-1}$.

We call a sequence $\Z$, with the above described partition and
associated sequence space $X^p$, a \term{general interpolating
sequence for $H^p$} if for each sequence $(w_k) \in X^p$ there is a
function $f\in H^p$ that interpolates it (in the sense described above).

If we choose a different norm than the sup norm for $\| \cdot \|_{G_k}$
then the space $X^p$ would not be the same (in general) and it would be
conceivable that the collection of general interpolating sequences
would change. However, it turns out that same collection results for any
reasonable choice of norms. This is explored a little more fully in
\cite{Lue04b}. Here, suffice it to say that the facts outlined in the
rest of this section remain valid if we changed, for example, to
$L^p(dA)$ means over $G_k$ or $L^p(ds)$ means over the boundary of
$(\Z_k)_\eps$. This is not surprising, since these norms are
dominated by the sup norm and also dominate the sup norm over slightly
smaller sets.

We list a number of facts about general interpolating sequences. We
omit the proofs where they are essentially the same as those in
\cite{Lue04b}.

\begin{fact}
  General interpolating sequences for $H^p$ are zero sequences for
  $H^p$.
\end{fact}

\begin{fact}\label{separation}
  If $\Z = \Union \Z_k$ is a general interpolating sequence for
  $H^p$ then there exists $\delta > 0$ such that for all $j \ne k$,
  $\psi(\Z_j, \Z_k) > \delta$.
\end{fact}

That is, while we don't necessarily have uniform discreteness for $\Z$
and one can have $\psi(z_{k_i},z_{j_i}) \to 0$, it can only happen if
eventually $z_{k_i}$ and $z_{j_i}$ are in the same $\Z_{m_i}$.

\begin{fact}\label{subsequence}
  A subsequence of a general interpolating sequence is also a
  general interpolating sequence. In particular, selecting one point
  from each $\Z_k$  produces an ordinary interpolating sequence.
\end{fact}

A consequence of this (and the requirements on the diameters of $G_k$)
is that there is an upper bound on the Carleson norms of all measures of the
form $\sum (1 - |a_k|)\delta_{a_k}$, where $a_k \in G_k$ is arbitrary.
Therefore $\sum_k \| f|_{G_k} \|_\infty d_k  \le C \| f \|_{H^p}$ and so:

\begin{fact}
  For a general interpolating sequence, the interpolation operator
  $T$ that takes $f \in H^p$ to the sequence $(w_k)$ where $w_k$ is the
  equivalence class of $f|_{G_k}$, is bounded onto $X^p$.
\end{fact}

If the interpolation operator is onto and bounded, then by the open
mapping theorem there exists an interpolation constant $M$. That is, if
$w \in X^p$ then there is a function $f_w\in H^p$ with norm
at most $M\| w \|$ that interpolates it.

\begin{fact}\label{expand}
  If $\Z = \{ z_1, z_2, \dots \}$ is a general interpolating
  sequence with constant $M$, $\eta > 0$, and $z_0$ is a point with
  $\psi(z_0, \Z) > \eta$, then $\Z \union \{ z_0 \}$ is also a
  general interpolating sequence \parenth{selecting any disk
  $D(z_0,r)$ for the corresponding $G_0$} with interpolation
  constant depending only on $M$, $p$, $\eta$ and $r$.
\end{fact}

\begin{fact}\label{upper}
  For any general interpolating sequence $\Z$, there is an upper
  bound $N$ on the cardinality of every $\Z_k$.
\end{fact}

This last fact follows from the Moebius invariance of interpolating
sequences. If it were not true we could take a $\Z_k$ with large
cardinality, and use a Moebius transformation to move it near the
origin. Then, using Fact~\ref{expand} if necessary, we could
obtain a function which is $1$ at $0$, zero at all points of
$\Z\setminus\{ 0 \}$ and with $H^p$ norm less than the interpolation
constant of $\Z \union \{ 0 \}$. This leads to a contradiction if the
cardinalities are unbounded

\section{The remaining proof}

From Facts~\ref{subsequence} and \ref{upper}, it follows that a
general interpolating sequence is a union of a finite number of
ordinary interpolating sequences. Note that this necessary condition is
independent of the choices involved: choosing how to divide $\Z$ into
$\Z_k$ and choosing what open sets $G_k$ are used.

What remains to be shown is the converse, that a Carleson sequence is a
general interpolating sequence.

A Carleson sequence automatically has bounded density (that is, there is
an upper bound on the number of points in any disk $D(z,R)$ for fixed
$R$) and so, by \cite{Lue04b}, there exists a partition of $\Z$ into
finite sets $\Z_k$ with these requirements satisfied: the sets $G_k =
(\Z_k)_\eps$ are disjoint for some $\eps > 0$ and $\diam_\psi(G_k) \le
C$ for some $C$ independent of $k$. Thus, the arguments of the rest of
this section will not be vacuous.

The $H^\infty$ case has essentially been know from the beginning. Any
sequence $(w_k) \in X^\infty$ corresponds to a bounded analytic function
on $G = \bigcup_k G_k$. Arc length on the boundary of $G$ is easily seen
to be a Carleson measure. It is then routine to show that $(w_k)$ can be
interpolated. The following merely fills in some of the details

Let $\Z_k$ and $G_k$ be as described. Let $w$ be an element of $X^p$ and
for each $k$ let $g_k$ be an analytic function on $G_k$ in the
equivalence class of $w_k$, such that $\| g_k \|_{G_k} = \| w_k \|_k$.
Clearly there exists, for any $J$, an element $p_J$ of $H^\infty$ which
interpolates $w_k$ for $1 \le k \le J$ and which is in fact a
polynomial. Using the usual duality methods, we can estimate the minimal
$H^\infty$ norm of such interpolators by
\begin{equation*}
  \sup \left| \int_{\IT} h(z)p_J(z)\overline{B(z)} \,dz \right|
\end{equation*}
where $B$ is the finite Blaschke product with zero sequence
$\Union_{k=1}^J \Z_k$, and the supremum is taken over all functions $h$
in the unit ball of $H^1$. Shrinking $\eps$ if necessary, we may
assume the sets $(\Z_k)_\eps$ are disjoint (Fact~\ref{separation}).
Let $\Gamma_k$ be the boundary of $(\Z_k)_\eps$. By Cauchy's
theorem, the above integral can be rewritten as the following integral
over the set $\Union_{k=1}^J \Gamma_k$, made up of of piecewise smooth
closed curves.
\begin{equation*}
  \sup \left| \int_{\IT} h(z)p_J(z)\overline{B(z)} \,dz \right| = \sup
    \left| \sum_{k=1}^J \int_{\Gamma_k} \frac{h(z)p_J(z)}{B(z)} \,dz
    \right|
\end{equation*}
Since the integrals in the summand depend only on the equivalence classes of
$(hp_J)|_{G_k}$ we get
\begin{align*}
  \sup \left|\int_{\IT} h(z)p_J(z)\overline{B(z)} \,dz \right|
    &\le \sup \left| \sum_{k=1}^J\int_{\Gamma_k} \frac{h(z)g_k(z)}{B(z)}
         \,dz \right|
    \\
    &\le \sup_{1\le k\le J} \frac{\| w_k \|}{\inf_{z\in \Gamma_k}
      |B(z)|} \int_{\Union\Gamma_k} |h(z)| \,dz.
\end{align*}
The infimum in the denominator above is bounded away from $0$
independent of $J$ for any ordinary interpolating sequence. Since $\Z$
is a finite union of ordinary interpolating sequences, that is true here
as well. Let $\delta > 0$ be that lower bound. Moreover, it is clear
that arc length on $\Union\Gamma_k$ is a Carleson measure with a constant
bounded above independent of $J$. Let $C < \infty$ be that upper bound.
Therefore, the minimal norm $f_J \in H^\infty$ that interpolates $w_k$
for $1 \le k \le J$ satisfies
\begin{equation*}
  \| f_J \|_\infty \le \frac{C}{\delta} \sup_k \| w_k \|.
\end{equation*}
Taking a limit point as $J \to \infty$ gives us an element of $H^\infty$
that interpolates all $w_k$.

This proves the $p = \infty$ case in almost exactly the same manner as
for ordinary interpolating sequences.

Now we turn to the case $0 < p < \infty$. For ordinary interpolating
sequences one can construct an interpolation operator from a bounded
sequence of bounded analytic functions $F_k$ that satisfy  $F_k(z_k) = 1$
and $F_k(z_j) = 0$ for $j \ne k$. This sequences is easily provided by the
Blaschke product: $F_k(z)=c_kB_\Z/\phi_{z_k}$ for appropriate constant
$c_k$. In our case we need something similar: $F_k|_{G_k}$ is in an
appropriate equivalence class, while $F_k|_{G_j}$ is equivalent to $0$
for $j \ne k$. This is easily provided by the solution of the $H^\infty$
case.

Let $w = (w_k) \in X^p$, so that $\| w \|_{X^p} = \sum \| w_k \|_k^p d_k
< \infty$. For each $k$ choose some $a_k \in G_k$, and consider the sum
\begin{equation}\label{VGHsum}
    f(z) = \sum_k  F_k(z) \left( \frac{1 - |a_k|^2}{1 - \bar a_k z}
                         \right)^2 e^{\beta_k(a_k) - \beta_k(z)}
\end{equation}
where $F_k \in H^\infty$ will be chosen later and where
\begin{equation*}
  \beta_k(z) = \sum_{j\ge k} (1 - |a_j|^2)\frac{1 + \bar a_j z}{1 - \bar a_j z}.
\end{equation*}
The summand in \eqref{VGHsum}, apart from the $F_k$, is essentially that
of Vinogradov, Gorin, and Hru\v s\v c\"ev in \cite{VGH81} where it was
shown that there exists a constant $C$, depending only on the Carleson
norm of the measure $\nu = \sum (1 - |a_k|) \delta_{a_k}$, such that
\begin{equation}\label{L-inftyEstimate}
  \sum_k \left| \frac{1 - |a_k|^2}{1 - \bar a_k z} \right|^2
        e^{\re(\beta_k (a_k) - \beta_k(z))} < C
\end{equation}
It was also shown that $\re\beta_k(a_k)$ is bounded, so the
exponent $\re(\beta_k (a_k) - \beta_k(z))$ in the above expression
is bounded above uniformly in $k$. It is a consequence of Harnack's
inequality that it is also bounded below when restricted to $G_k$,
uniformly in $k$. It is clear that
\begin{equation}\label{L-oneEstimate}
  \frac{1}{2\pi} \int_0^{2\pi} \frac{1 - |a_k|^2} {|1 - \bar a_k
    re^{i\theta}|^2} \le 1.
\end{equation}
Setting $z = re^{i\theta}$ we view
\begin{equation*}
  \frac{1 - |a_k|^2}{(1 - \bar a_k re^{i\theta})^2} e^{\beta_k (a_k) -
    \beta_k(re^{i\theta})}
\end{equation*}
as the kernel $G_r(\theta, \zeta)$ of an operator from $L^p(\nu)$ to
$L^p(d\theta)$ (recall that $\nu = \sum (1 - |a_k|^2)\delta_k$). It
satisfies (by \eqref{L-inftyEstimate} and \eqref{L-oneEstimate})
\begin{equation*}
    \int |G_r(\theta, \zeta)| \,d\nu(\zeta) < C \quad\mbox{and}\quad
    \int |G_r(\theta, \zeta)| \,d\theta < C,
\end{equation*}
and so it defines a bounded operator between $L^p(d\nu)$ and
$L^p(d\theta)$ for all $1 \le p \le \infty$. Thus we have
\begin{equation}\label{basicinequality}
  \int_0^{2\pi} \left( \sum_k  b_k \left| \frac{1 - |a_k|^2}{1 - \bar
            a_k re^{i\theta}} \right|^2 e^{\re[\beta_k(a_k)-
            \beta_k(re^{i\theta})]}
            \right)^p \,d\theta \le C \sum_k b_k^p d_k
\end{equation}
for all sequences $b_k \ge 0$.

Give $(w_k) \in X^p$, we choose $F_k$ so that $F_k|_{G_j}$ is equivalent
to zero when $j\ne k$ and such that
\begin{equation*}
  F_k(z) \left( \frac{1 - |a_k|^2}{1 - \bar a_k z} \right)^2
            e^{\beta_k(a_k)- \beta_k(z)}
\end{equation*}
is equivalent to a $g_k$ in the class $w_k$. We can do this by choosing
$F_k(z)$ equivalent to
\begin{equation*}
  g_k(z) \left/ \left( \left( \frac{1 - |a_k|^2}{1 - \bar a_k z} \right)^2
    e^{\beta_k(a_k)- \beta_k(z)}\right) \right.
\end{equation*}
Since the denominator above has absolute value bounded below in $G_k$
(independent of $k$) this can be done with $\| F_k \|_\infty < C \| w_k
\|_k$ by the $p = \infty$ case. Now, since
\begin{align*}
   |f(z)| &= \left| \sum_k F_k(z) \left( \frac{1 - |a_k|^2}{1 - \bar a_k
                 z} \right)^2 e^{\beta_k(a_k) - \beta_k(z)} \right|\\
          &\le C \sum_k  \| w_k \|_k \left| \frac{1 - |a_k|^2}{1 - \bar
                a_k z} \right|^2 e^{\re[\beta_k(a_k)- \beta_k(z)]}
\end{align*}
we have $\| f \|_{H^p} \le \| w \|_{X^p}$ by the inequality
\eqref{basicinequality}. Clearly also $f$ interpolates $w$ by the
construction of $F_k$.

It remains to consider the case $p < 1$. In this case we perform the
same construction with the sum
\begin{equation*}
  f(z) = \sum_k  F_k(z) \left( \frac{1 - |a_k|^2} {1 - \bar a_k z}
            \right)^{2/p} e^{[\beta_k(a_k)- \beta_k(z)]/p}
\end{equation*}
Then estimating $\| f \|_{H^p}$ is a simple matter of
taking all terms in the sum to the power $p$ and integrating term by
term. In all cases we obtain a function $f \in H^p$ that interpolates
the given $w\in X^p$, and so $\Z$ is a general interpolating
sequence. This completes the proof.

%

\newcommand{\noopsort}[1]{}
\providecommand{\bysame}{\leavevmode\hbox to3em{\hrulefill}\thinspace}
\providecommand{\MR}{\relax\ifhmode\unskip\space\fi MR }
\providecommand{\MRhref}[2]{%
  \href{http://www.ams.org/mathscinet-getitem?mr=#1}{#2}
}
\providecommand{\href}[2]{#2}

\end{document}